\documentclass[12pt]{amsart}
\usepackage[utf8]{inputenc}
\usepackage{graphicx}
\usepackage{tikz}        
\usepackage{amssymb}
\usepackage{epstopdf}
\usepackage[mathscr]{eucal}
\usepackage{amsmath,amssymb,amscd}
\usepackage{epsfig}
\usepackage{bbold}
\usepackage{hyperref}
\usepackage{MnSymbol}
\usepackage{amsfonts}
\usepackage{xcolor}%

\usepackage[numbers, sort&compress]{natbib}

\newtheorem{theorem}{Theorem}
\newtheorem{lemma}{Lemma}

\newtheorem{proposition}{Proposition}
\theoremstyle{definition}
\newtheorem{remark}{Remark}

\def \mb{\mathbb}

\def \mc{\mathcal}

\def \Z{\mb Z}                  
\def \R{\mb R}                 

\def \a{\alpha}         
\def \D{\Delta}         
\def \vp{\varphi}       
\def \th{\theta}       

\def \S{\mb S}        

\title{The Schubart orbits on the circle }

\begin{document}
	\maketitle
	\markboth{Shuqiang Zhu}{  Schubart orbits on $
		\S^1$}
	\vspace{-0.5cm}
	\author       
	\bigskip
	\begin{center}
		{Shuqiang Zhu}\\
		{\footnotesize 
			
			School of  Mathematics,  Southwestern University of Finance and Economics, \\
			
			Chengdu 611130, China\\
			
			zhusq@swufe.edu.cn 
			
		}
	
	\end{center}

\begin{abstract}
We consider the three body problem on $\S^1$ under the cotangent potential.  We first construct homothetic orbits ending in singularities, including total collision singularity and collision-antipodal singularity. Then certain   symmetrical periodic orbits with two equal masses, called Schubart orbits, are shown to exist. The proof is based on  the construction of a Wazewski set  in the phase space. 
\end{abstract}

\section{Introduction}

The Newtonian n-body problem has been generalized in many ways, for instance,  under the general homogeneous potentials, or  in higher dimensional Euclidean spaces.  Among them, the curved n-body
problem,  which studies  n-body problem on  surfaces of constant curvature under the cotangent potential, 
has received lot of attentions in the last decade (cf \cite{borisov2018reduction, diacumemoir, diacu2018central} and the references therein ).

In the Newtonian n-body problem, the two-body case is a Hamiltonian system with one degree of freedom, so is integrable.  The cases with  two degrees of freedom, namely, the restricted three-body problem and the collinear three-body problem, remains to be largely unsolved.
However, the ideas emerged in attacking them shed light on more general problems (cf \cite{Moser1973, Mcgehee1974, Xia5body} and the references therein).   

We consider   one case  of the curved n-body problem with   two degrees of freedom, namely,  the three-body problem on the circle.  As a  preliminary study on the problem, in this manuscript we construct some interesting orbits, the Eulerian homothetic orbits and Schubart orbits.  

The Eulerian homothetic orbits are connected with the singularities.  In the curved three-body problem, there are collision, antipodal and collision-antipodal  singularities. The antipodal singularities   turn out to be impossible by consideration of Hill's region.  The collision-antipodal  singularities are more interesting. In fact,  our study indicates that it might be sensitive to the choice of masses. 

Intuitively, if all the masses stay on a semi-circle, they would attract each other and end in a triple collision. However, if two masses are equal, there exist certain symmetric periodic orbits,  provided collisions are regularized. The behavior of   these orbits is as follows. Initially,  the unequal mass, $m_3$, is at the midpoint of the equal ones.  If the masses were released with zero velocity, it would be a homothetic collapse to triple collision.   However, we set the initial velocities such that  $m_2$ and $m_3$ move towards each other and   $m_1$ leaves them. Then $m_1$ slows and stops exactly when $m_2$ and $m_3$ collide. This is the first quarter of the orbit.  The second quarter of the orbit is the time-reverse of the first, and the second half is the reflection of the first half with the roles of $m_1$ and  $m_2$ reversed.   Such orbits are called \emph{Schubart orbits}. They were found numerically by Schubart \cite{Schubart1956} for the Newtonian three-body problem and the analytic existence proof, was given by in  \cite{Moeckel-Schubart,Shibayama2011}, among others.  

For the  existence proof  of the Schubart orbits, we follow that of Moeckel in \cite{Moeckel-Schubart}. It is a topological argument and it is  a variation of an idea used by Conley \cite{Conely_retrograde} in the Newtonian restricted three-body problem.  More precisely, it is based on the construction of a Wazewski set. 
Unlike that of the Newtonian case, where the potential is almost a function of the shape variable, the cotangent potential depends essentially on the two variables. Some computations are relatively lengthy.

The paper is organized as follows. In Sect. \ref{sec:setting}, we discuss the basic setting of the three-body problem on $\S^1$ and the Eulerian homothetic orbits. In Sect. \ref{sec:regular}, we regularize the collision singularities.  In Sect. \ref{sec:shooting}, we apply a topological argument to show the existence of Schubart orbits. Some technical computations are presented   in the Appendix.

\section{Settings and Eulerian homothetic orbits}\label{sec:setting}
The configuration space is  $(\mathbb S^1)^3$.  The coordinates are $\vp_1, \vp_2, \vp_3$, with  $\varphi_i\in \R/2\pi \Z. $
The Lagrangian is 
\[L=\sum_i \frac{1}{2} m_i \dot {\varphi}^2_i+\sum_{j\neq i} m_im_j\cot d_{ij}, \]
where  $d_{ij}=\min \{|\vp_i -\vp_j|, 2\pi -|\vp_i -\vp_j| \}$.  The system is undefined in the set $\D=\bigcup_{i\ne j} \D_{ij}$, with 
\[  \D_{ij}=\{ (\vp_1, \vp_2, \vp_3): d_{ij}=0   \} \bigcup \{ (\vp_1, \vp_2, \vp_3): d_{ij}=\pi   \}.   \]
The cases $d_{ij}=0$ are  collisions, whereas the cases $d_{ij}=\pi$ are   antipodal configurations, when some bodies are at the opposite ends of a diameter. In both cases, the forces are  infinite. There are other possibilities.  For instance,  consider the case $d_{12}=0, d_{23}=\pi$, which corresponds to a configuration with $m_1, m_2$ at collision and $m_3$ lies at the opposite end. This  will be called \emph{collision-antipodal singularity}, \cite{diacu2011, diacu2012n2}.

The rotation group $SO(2)$ acts on the configuration space by 
$$(\vp_1, \vp_2, \vp_3) \mapsto (\vp_1+s, \vp_2+s, \vp_3+s),   s\in \R. $$ This action keeps the potential function, and actually keeps the system by  tangent lift. 
The corresponding first integral  is the angular momentum, i.e.,  $\mathbf J= \sum_{i=1}^n m_i \dot\vp_i$. 
Thus, we have two first integrals
\begin{equation}\label{equ:firstintegral}
\sum m_i \vp_i= \a t + \a', \ \ K -U =h.   
\end{equation}
Note that  is $\vp_i(t)$ is a solution, so is $\vp_i(t) + a t + a'$. 
We may assume, by changing the coordinates by linear functions of the variable $t$, that  
\[   \sum m_i \vp_i=0 \ (\mod 2\pi). \]
We further  assume that $\sum m_i=1$. 

\subsection{Jacobi Coordinates and the Singularities}

Let 
\[ 1=m_1+m_2+m_3, \ \ \ \ \a_1= \frac{m_1}{m_1+m_2}, \ \ \ \a_2= \frac{m_2}{m_1+m_2}.\]
Introduce Jacobi variables as 
\[ x_1=\varphi_2-\varphi_1, \ \ \ x_2=\varphi_3-\a_1\varphi_1-\a_2\varphi_2  \]
and their velocities $u_i =\dot x_i$. Then 
the inverse is 
\begin{align*}
& \vp_1 =- \a_2 x_1 - m_3x_2, \ & \vp_2 =\a_1 x_1 - m_3 x_2,  \ \  \ & \vp_3 =(m_1+m_2)x_2, \\
& \vp_2-\vp_1= x_1,\   &\vp_3-\vp_1= \alpha_2 x_1+ x_2,\  \ \  &\vp_2-\vp_3= \alpha_1 x_1- x_2.  
\end{align*}

Then the kinetic energy is 
\begin{align*}
2K &= \sum m_i \dot \vp_i^2 = m_1(- \a_2\dot  x_1 - m_3 \dot x_2 )^2 + m_2 ( \a_1 \dot x_1 -  m_3 \dot x_2)^2 + m_3 (m_1+m_2)^2 \dot x_2 ^2\\
&=\mu_1 u_1^2 + \mu_2 u_2^2, 
\end{align*}
where $ \mu_1=\frac{m_1m_2}{m_1+m_2}$ and  $\mu_2=(m_1+m_2)m_3.$  The potential is
\[  U(x_1, x_2)=m_1 m_2 \cot d_{12}  +m_1 m_3 \cot d_{13}  +m_2 m_3 \cot d_{23} .   \]

Now we assume that the three bodies are ordered on the circle anti-clockwise as 
$$-\pi+2k\pi \le \vp_1\le  \vp_3\le  \vp_2\le \pi + 2k \pi. $$
Then
\begin{align*}
&d_{12}=\min\{ x_1, 2\pi-x_1\}, \   d_{13}=\min\{ x_2+ \a_2 x_1, 2\pi-x_2- \a_2 x_1\}, \\
&  d_{23}=\min\{ \a_1 x_1-x_2, 2\pi-\a_1 x_1+x_2\}. 
\end{align*}

 For the Newtonian collinear three-body problem,  a similar ordering is $q_1<q_3<q_2$, where $q_i$ is the coordinates of $m_i,$ $i=1, 2, 3$. 
For such an ordering,   possible singularities are 
\emph{total collision}, \textit{collision between} $m_2, m_3$,  \textit{collision between} $m_1, m_3$.   With the Jacobi coordinates, the configuration space is a region between  two half lines, \cite{Moeckel-Schubart}.

\begin{figure}[!h]
	\centering
	\includegraphics[scale=0.5]{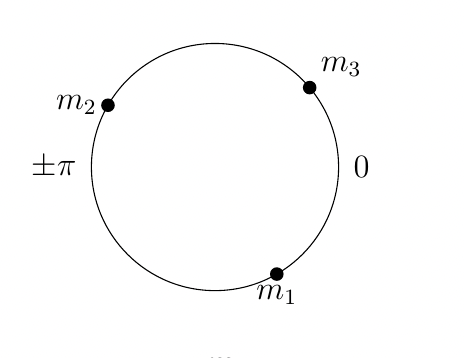}
	\includegraphics[scale=0.6]{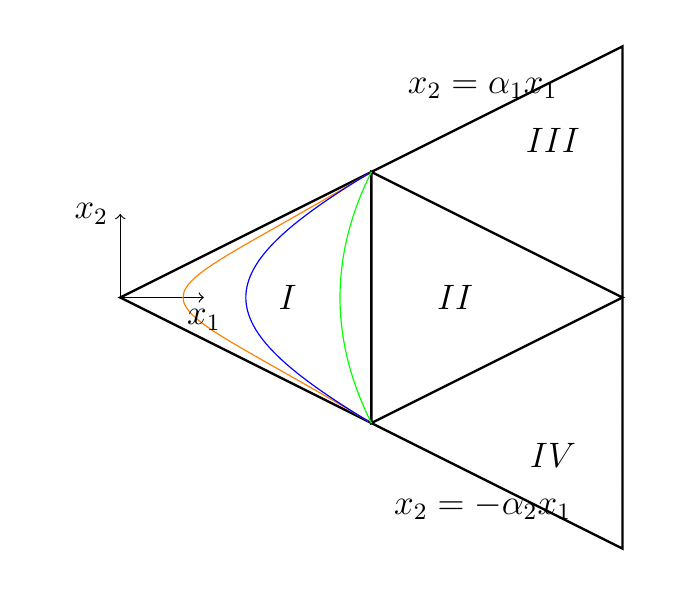}
	
	\caption{The configuration space, and some zero velocity curves. The masses are  $\frac{1}{3}, \frac{1}{3}, \frac{1}{3}$. Zero velocity curves in region I for $h= -100$ (orange), $0$ (blue), $100$ (green)  are shown.}
	\label{fig:region-s1}
\end{figure}

For  our  three-body problem on $\S^1$,  obviously, we have 
\[   \vp_2-\vp_1 = x_1\in [0, 2\pi], \ \vp_2-\vp_3= \a_1 x_1 -x_2\in [0, 2\pi], \ \ \vp_3-\vp_1= \a_2 x_1 +x_2\in [0, 2\pi],     \]
Then  the configuration space is  the triangular region bounded by 
\[ \a_1 x_1 -x_2=0, \ \a_2 x_1 +x_2=0, \  x_1 = 2\pi.    \]
as shown in Figure \ref{fig:region-s1}. 
The singularities are 
\[ \vp_2-\vp_1, \vp_2-\vp_3 , \vp_3-\vp_1=0, \pi, 2\pi.     \]
That is,  
\[  x_1 = 0, \pi, 2\pi, \ \a_1 x_1 -x_2= 0, \pi, 2\pi,  \ \a_2 x_1 +x_2= 0, \pi, 2\pi, \]
i.e., the vertices, boundary and the three mid-segments of the triangular region.  These singularities  divide the configuration space into four sub-triangles.   Let us denote  the four sub-triangular regions as  I, II, III, and IV, as in Figure \ref{fig:region-s1}. 

	

\begin{figure}[h!]
	\centering
	\includegraphics[scale=0.3]{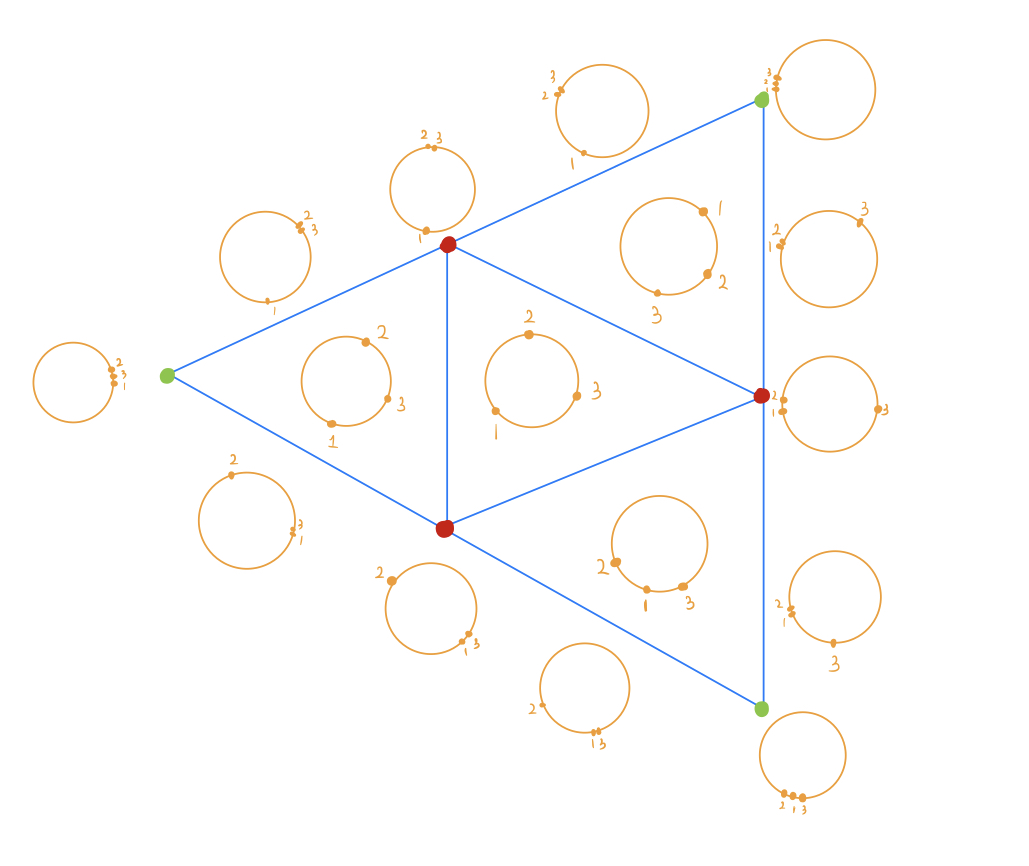}
	
	\caption{The configurations corresponding to vertices, sides and interior of the four sub-triangular regions }
	\label{fig:configuration}
\end{figure}

The three vertices correspond to the total collisions, the three sides correspond to the double collisions, the three mid-segments correspond to 
antipodal singularities, and the intersections of the three mid-segments correspond to 
three collision-antipodal singularities. 
 In Figure \ref{fig:configuration}, we sketch the real configurations corresponding to typical points of the configuration space.

\subsection{Hill's region}
Consider the motion on the energy surface $H=h$, then the projection of the energy surface  to the configuration space is called the Hill's region corresponding to energy $h$, 
\[  \mathcal H(h)= \{ x: U(x)+h\ge 0  \} \]
Recall that the energy surface  lies over its projection $\mathcal H(h)$ as a kind of degenerate circle bundle. The boundary  $\partial \mathcal H(h) = \{x: U(x) = -h\}$ is  the zero-velocity curve.

We sketched some zero-velocity curves in region I in Figure \ref{fig:region-s1}. The zero-velocity curves in region  III and IV are similar to that in region I, but they are more complex  in region II, as we shall see soon.     Note that $U$ is undefined at the three intersections of the mid-segments. Other than the three points, $U(x)\to \infty$ as $x$ approaches the three sides, and  $U(x)\to -\infty$ as $x$ approaches the three mid-segments. 
Thus we have the following 

\begin{proposition}
The  antipodal singularities  are repelling for the three-body problem on $\S^1$. 
\end{proposition}

\subsection{Eulerian homothetic orbits}\label{subsec:homothetic}

For three-body problem, it is enough to consider the motion in region I and II.  Let us first consider a simple case.

 Consider the isosceles problem on $\S^1$. The masses are $m_1=m_2=\frac{1-m}{2}, m_3=m$. The  initial data is 
 \[ \varphi_2=\varphi \in (0, \pi),   \varphi_1= -\varphi,  \varphi_3=0,   \dot \varphi_1= -\dot \varphi_2, \dot \varphi_3= 0.   \] 
By the symmetry, the configuration would stay isosceles with $\vp_3=0$.
 That is, the system has just one degree of freedom.  More precisely, since $\vp_2-\vp_3=\vp_3-\vp_1$, we see $x_2=0$, and $x_1=2\vp \in (0, 2\pi)$, $\mu_1=m_1/2=\frac{1-m}{4}$. Let $p_1=\mu_1u_1$, then
\[ \ H =\frac{1}{2} \frac{p_1^2}{\mu_1} - U(x_1),  \  U=\frac{(1-m)^2}{4}  (\cot  d_{12} + \frac{4m}{1-m}  \cot d_{23}).     \]
The motion depends on the   function, 
\[  \begin{cases}
  \cot  x_1 + \frac{4m}{1-m}  \cot \frac{x_1}{2},    \  \ \ 0< x_1<\pi \cr
- \cot  x_1 + \frac{4m}{1-m}  \cot \frac{x_1}{2}.    \  \ \ \pi < x_1<2\pi, \cr
\end{cases}\]
On $[0,\pi]$, the function is decreasing from $\infty$ to $-\infty$ for any value of $m$. While on $[\pi,2\pi]$,  the graph  depends  on the value of $m$, since 
 \begin{align*}
&- \cot 2\vp  +\frac{4m}{1-m}  \cot \varphi = \frac{\sin^2 \varphi-  \cos^2 \varphi }{ 2 \sin \varphi \cos \varphi  }  +\frac{4m}{1-m}  \cot \varphi\\
&=\frac{\sin^2 \varphi }{ 2 \sin \varphi \cos \varphi  }  +(\frac{4m}{1-m}  -\frac{1}{2}) \cot \varphi= \frac{1}{2}\tan \frac{x_1}{2} + \frac{9m-1}{2(1-m)} \cot \frac{x_1}{2}. 
\end{align*} 

  \begin{figure}[h!]
 	\includegraphics[scale=0.4]{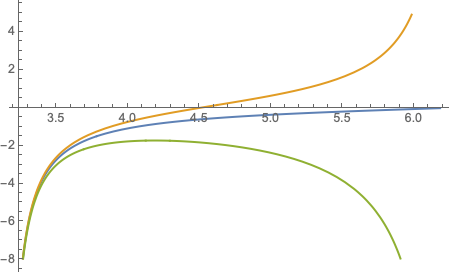}
 	\caption{ Graphs of $U$ on $[\pi,2\pi]$ for  $m<1/9$ (orange),  $m=1/9$ (blue),  $m>1/9$ (green)}
 	\label{fig:isocseles_2}
 \end{figure}

The motions are obtained by the conservation of energy. Recall  that $x_1=0,	\pi, 2\pi$ are   singularities. More precisely, $x_1=0$ is the total collision, 	$x_1=\pi$ is the antipodal singularity, and $x_1=2\pi$ is the collision-antipodal   singularity. 
On $(0, \pi)$, all motions would eventually go to $x_1=0$, or, triple collision must happens. At $x_1=\pi$, $-U$ is $\infty$, so it is repelling, or, antipodal singularity is impossible.  The qualitative feature of  motions on $(\pi, 2\pi)$ depends on the sign of 
$9m-1$. 
\begin{itemize}
	\item If $m<1/9$,  all motions would eventually go to $x_1=2\pi$, or, collision-antipodal singularity  must happen, and at that moment, the velocity is infinite.
	
	\item If $m=1/9$,  $-U$ is decreasing on $(\pi, 2\pi)$ and $-U(2\pi) =0$. Then for   $H=h\ge 0$,  collision-antipodal singularity  must happen, and at that moment, the velocity is finite.

	\item If $m>1/9$,  there is one critical value $h_0$ of $-U$.  If $H=h>h_0$, the motions are periodic. If $H=h=h_0$, it is a stable equilibrium. 
\end{itemize}

For later use,  let us refer to those configurations as \emph{Eulerian central configurations} and the orbits  as \emph{Eulerian homothetic orbits}. 
 \begin{remark}
 	The above example was first  considered by 
Florin et al, see \cite{diacu2012n2}. 
 \end{remark}

\section{Regularization of the collision singularities}\label{sec:regular}
We now focus on motions in  region I. Intuitively, all motions in this region seems to end in a total collision. However, we will construct symmetric periodic orbits in region I, called Schubart orbits, in next section.  In this section, we regularize the double and triple collision singularities.

Assume that $m_3=m\in (0,1)$, $m_1=m_2=n$, then $n=(1-m)/2$, and 
\[\alpha_1 =\alpha_2=1/2, \  \mu_1 = \frac{1-m}{4}= n/2, \ \mu_2=(1-m)m=2n m.  \]
In region I, the distances are $d_{12}=x_1, d_{13}=\frac{1}{2} x_1 +x_2, d_{23}=\frac{1}{2} x_1 -x_2$, so 
\begin{align*}
&L = K+U, \ K=\frac{1}{2}(\mu_1 u_1^2 + \mu_2 u_2^2),\\ 
&U(x_1, x_2)=n^2 \cot x_1  +m n \left(\cot(  \frac{1}{2} x_1 +x_2)  +  \cot (  \frac{1}{2} x_1-x_2) \right) 
\end{align*}

Recall that region I is a triangular region  bounded by $ \frac{1}{2} x_1 \pm x_2=0$ and $x_1=\pi$. The three vertices and sides are singularities.    We perform first  Mcgehee's coordinates then another change of variables to eliminate the singularities corresponding to the  collisions, see Figure \ref{fig:region-s1}. 


Let 
\[  x_1 =\frac{1}{\sqrt{\mu_1}} r \cos \th,  \ x_2 =\frac{1}{\sqrt{\mu_2}} r \sin \th.   \]
Then $2K=\dot r^2 + r^2 \dot \th^2$, and 
\[ U(x_1, x_2)=n^2 \cot (rA_1 \cos \th )  +m n \left(\cot[ r  A_2 \sin(\th+\th_*)] +  \cot[ r  A_2 \sin(\th_*-\th)] \right)
\]
where 
$  A_1= \sqrt{\frac{ 2}{n}}, A_2 =\sqrt{\frac{m+1}{2 n m}}, $ and 
\[   \th_* = \arctan\sqrt{m}, \  \th_* < \frac{\pi}{4}  \]
The configuration space has been blew up to
\[  \th \in (-\th_*, \th_*), \  0\le r < \frac{\sqrt{\mu_1}\pi}{\cos \th}= \frac{\sqrt{n }\pi}{\sqrt{2}\cos \th}.   \]

The corresponding second-order Euler-Lagrange equations are:
\begin{equation}\label{equ:mcgehee_1}
\begin{split}
\ddot r&=r \dot \th^2 + U_r, \\
\dot{\overline{ r^2 \dot \th}}&= U_\th. 
\end{split}
\end{equation}
Next, one can blow-up the triple collision singularity at $r=0$ by introducing the 
the time rescaling $'=r^{\frac{3}{2}}\  \dot {}$ and the variable $\nu  = r'/r$. Setting $\tau =\th'$  gives the following first-order system of differential equations:
\begin{equation}\label{equ:mcgehee_2}
\begin{split}
r' &= \nu r, \\
\nu'&=\frac{1}{2} \nu^2 + \tau^2 + r^2 U_r\\
\th '&= \tau \\
\tau ' &=-\frac{1}{2} \tau \nu + r U_\th
\end{split}
\end{equation}
with energy equation:
\begin{equation}\label{equ:energy_mcgehee_1}
\frac{\nu^2 +\tau^2}{2}- rU= rh
\end{equation}
Explicitly,  the functions are 
\begin{align*}
rU&= n^2 r\cot (rA_1 \cos \th )  +m n \left(r\cot[ r  A_2 \sin(\th+\th_*)] +   r\cot[ r  A_2 \sin(\th_*-\th)] \right)\\
rU_\th&=n^2\frac{r^2 A_1 \sin \th}{\sin^2 (rA_1 \cos \th )}  +m n \left(  -\frac{r^2  A_2 \cos (\th + \th_*)}{\sin^2[ r  A_2 \sin(\th+\th_*)]} +    \frac{r^2  A_2 \cos (\th_*- \th)}{\sin^2[ r  A_2 \sin(\th_*-\th)]} \right)\\
r^2 U_r&=-n^2\frac{r^2 A_1 \cos \th}{\sin^2 (rA_1 \cos \th )}   -m n \left(  \frac{r^2  A_2 \sin (\th + \th_*)}{\sin^2[ r  A_2 \sin(\th+\th_*)]} +    \frac{r^2  A_2 \sin (\th_*- \th)}{\sin^2[ r  A_2 \sin(\th_*-\th)]} \right)
\end{align*}
They are well-defined at $r=0$. Hence $\{r = 0\}$ is now an invariant set for the flow, called the \emph{triple collision manifold}.

The differential equations are still singular due to the double collisions at $\th =\pm \th_*$.  The final coordinate change will eliminate these singularities.
Define new variables $u, \gamma$ such that
\[ \th = \th_*  \sin u, \  \gamma   = \tau \cos^2 u.  \]
Note that $-\th_* \le  \th \le \th_*$ corresponds to  $-\frac{\pi}{2} \le u \le \frac{\pi}{2}$. After calculating the differential equations for $u, \gamma$,  introduce a further rescaling of time by multiplying the vector field by $\th_* \cos^2 u$. Retaining the prime to denote differentiation with respect to the new time variable one finds 
\begin{equation}\label{equ:mcgehee_3} 
\begin{split}
r' &= \th_*\nu r\cos^2 u \\
\nu'&= \th_*\cos^2 u(\frac{1}{2} \nu^2 + \tau^2 + r^2 U_r)=\th_*(\frac{1}{2} \nu^2\cos^2 u + \frac{\gamma^2}{\cos^2 u} +  r^2 U_r \cos^2 u )\\
&=\th_*\cos^2 u(-\frac{1}{2} \nu^2 + 2 rh + 2rU + r^2 U_r)\\
u '&=\frac{ \gamma}{\cos u}\\
\gamma ' &=-\frac{1}{2}\th_* \nu \gamma \cos^2 u  + \th_* r U_\th \cos^4 u -2 \tau^2 \sin u\cos^2 u,   \\
&=-\frac{1}{2}\th_* \nu \gamma \cos^2 u  + \th_* r U_\th \cos^4 u -2  \sin u \frac{\gamma^2}{\cos^2 u}, 
\end{split}
\end{equation}
with energy equation:
\begin{equation}\label{equ:energy_mcgehee_2}
\frac{\nu^2\cos^2 u +\gamma^2/\cos^2 u}{2}- r U \cos^2 u= rh\cos^2 u
\end{equation}
  The configuration space is now 
\[  u \in \R, \ 0\le r < \frac{\sqrt{n }\pi}{\sqrt{2}\cos (\th_* \sin u)}. \]

Note that there is still one singularity on the boundary of the configuration space,  $r= \frac{\sqrt{n }\pi}{\sqrt{2}\cos (\th_*)}, u=\frac{\pi}{2}.$ Recall that it is one of the intersections of the mid-segments and that the potential is undefined there.  Denote   it by $Q$, see Figure \ref{fig:rectangle}.  Except this singularity,  The vector field is smooth and continuous on the boundary.

The differential equations \eqref{equ:mcgehee_3} represent the three-body problem on $\S^1$ with the prescribed energy for configurations being an obtuse triangle and  with $m_3$ in the middle, with triple collision blown-up and double collisions regularized. The shape variable $u$  need not be restricted to the interval $[-\frac{\pi}{2},\frac{\pi}{2}]$.  As $u$ ranges over the real axis,  the configuration oscillates between the double collisions at $\pm \th_*$ and the mass $m_3$ collides with $m_1$ and $m_2$ successively.

The equations  has some symmetries and 
a  Schubart orbit can be obtained from an orbit from $u=0,\nu=0$ to $u=\frac{\pi}{2}, \nu=0$. Suppose that we have a quarter of the trajectory $\Gamma(t), t\in [0,t_1]$ with 
\[  \Gamma(0)=(r_0, 0, 0, \gamma_0), \  \Gamma(1)=(r_1, 0, \frac{\pi}{2}, \gamma_1).  \] 
we can construct the second quarter of the  trajectory from $\Gamma(1)$
to $ \Gamma(2)=(r_0, 0, \pi, \gamma_0)$  by 
\[ [t_1, 2t_1] \to (r(2t_1-t), -\nu(2t_1-t), \pi-u(2t_1-t), \gamma(2t_1-t)),   \] 
(it satisfies the boundary condition and is a solution since  $ U_\th(\pi-u)=U_\th(u),  U_r(\pi-u)=U_r(u) $) and so the third and the fourth quarters. 
Then, using the symmetry of the vector field, it follows that one can piece together the four of them.  Hence, the existence of the required orbit is  reduced to find the first quarter, which  will be proved by a topological shooting argument in the next section.

\section{Schubart orbits by the shooting method}\label{sec:shooting}
Consider the system \eqref{equ:mcgehee_3} on the manifold of fixed energy $h = -1$.  We  construct the first quarter of the claimed Schubart orbit in this section. We will follow the shooting method in \cite{Moeckel-Schubart}, where Moeckel used it to show the existence of Schubart orbit for the Newtonian collinear three-body problem.  The idea is to construct a continuous map in the phase space and then apply a shooting argument. 

The construction of the  continuous map in the phase space is based on the result of Wazewski \cite{Wazewski1947}. Roughly speaking, a subset, called a \emph{Wazewski set},  of the phase space is carefully chosen such that the amount of time required to leave depends continuously on initial conditions. Then the exit point also depends continuously on initial conditions. This idea were developed by Conley and Easton \cite{Conley_isolating, Conly_book, Easton1970} to isolating blocks, topological index for invariant sets.

 There are several technical computations in this section. To not interrupt the flow of the argument, we will just claim them in this section and give the detail    in the appendix.

\subsection{The Wazewski set $\mc W$}
Consider a flow $\phi_t(x)$ on a topological space $X$ and a subset $\mc W\subset X$ of $X$.  Let   $\mc W_0$ be the set of points in $\mc W$ which eventually leave $\mc W$ in forward time, and let $\mc E$ the set of points which exit immediately:
\begin{align*}
&\mc W_0 =\{x\in \mc W: \exists t>0, \phi_t(x)\notin \mc W\}, \\
&\mc E=\{x\in \mc W: \forall t>0, \phi_{[0,t)}(x)\nsubseteq\mc W\}. 
\end{align*}
Clearly, $\mc E \subset \mc W_0$. Given  $x\in \mc W_0$ define the \emph{exit time}
\[  \tau (x) = \sup\{ t\ge 0: \phi_{[0,t)}(x)\subseteq\mc W\}\}.    \]
Note that $\tau(x) = 0$ if and only if $x\in \mc E.$ Then $\tau$ is continuous if 
\begin{itemize}
	\item    If $x\in \mc W_0$ and $\phi_{[0,t]}(x)\subseteq \overline{\mc W}$, then $\phi_{[0,t]}(x)\subseteq \mc W$.
	\item  $\mc E$ is a relatively closed subset of $\mc W_0$.
\end{itemize}
In this case, 
the set   $\mc W$ is called a \emph{Wazewski set}.

Since we are considering the motion on the energy manifold $h=-1$,  the configuration  is in the  region, $\{(r,u): U\ge 1\}$. 
 Note that $U=1$ defines an implicit function $r(u)$ since $U_r<0$. Define $r_*= r(0)$. 
Let
\[  \mc W = \{  (r, \nu, u,\gamma):\eqref{equ:energy_mcgehee_2} \ {\it holds}, \ 0\le r \le r_*, \ \nu \le 0, 0\le u\le \frac{\pi}{2}, \gamma \ge 0    \}, \]

\begin{figure}[h!]
	\centering
	\includegraphics[scale=0.2]{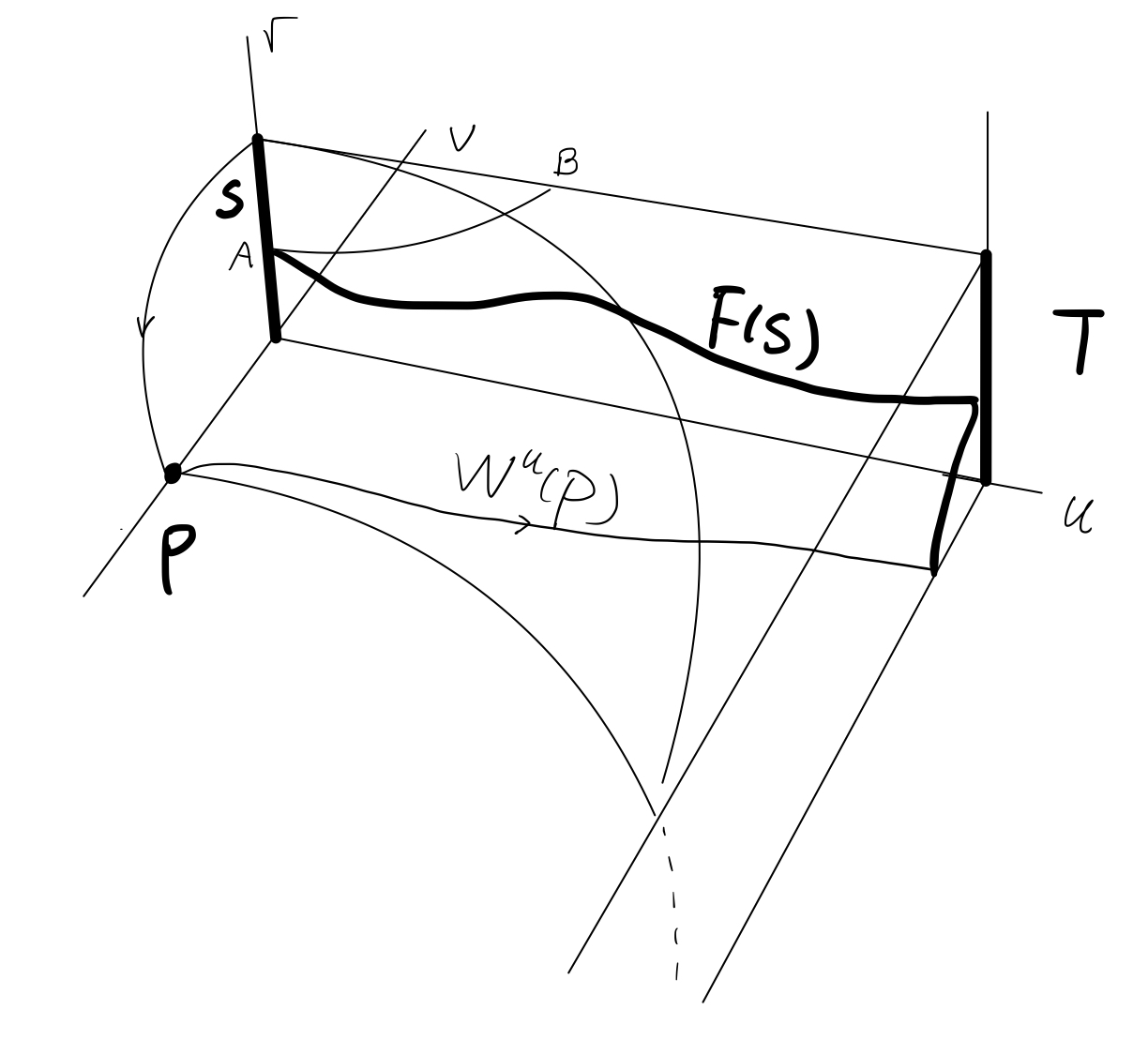}
	
	\caption{The Wazewski set $\mc W$ and a sketch of the existence proof }
	\label{fig:Waze-set}
\end{figure}

The choice  of the set $\mc W$ is motivated by that of Moeckel in \cite{Moeckel-Schubart}. The major  difference is that the value of $r$ is confined to $[0,r_*]$ in our case, and it is not in that of Moeckel's proof. Thus, the configuration space is restricted to a rectangle,  see Figure \ref{fig:Waze-set}. 
As it turns out, the restriction $0\le r \le r_*$ is essential for our proof. On one hand, this restriction   avoids  the 
singularity Q so that the system leads to a well-defined flow on $\mc W$. On the other hand,  the restriction leads to the  estimate\eqref{equ:estimate}, which is essential for our proof. Note that the restriction 
makes no harm since  $r$ is non-increasing in $\mc W$, so the exit points must have $r\le r_*$.

To visualize  $\mc W$,  we  use coordinates $(r, \nu, u)$ on the energy manifold, and the value of $\gamma$ is determined by
energy equation. The energy manifold projects to the three-dimensional  region
\[\frac{\nu^2\cos^2 u }{2}- r U \cos^2 u+ r\cos^2 u \le 0, \ 0\le r\le r_*,\]
see  Figure \ref{fig:Waze-set}. The south part of the  upper surface in the figure, where equality holds in \eqref{equ:energy_mcgehee_2} corresponds to $\gamma= 0$. The figure also  shows a sketch of the shooting argument.

\subsection{The invariant manifold $\mc H=\{ u=0, \gamma =0\}$. }

It is easy to verify that it is invariant under the flow since 
$ u'= \frac{\gamma}{\cos u}=0$ and 
\[  \gamma'= \th_* r U_\th|_{\th=0} \]
where we have used that fact that $\tau=0, \gamma=0$ and the following claim. 

\textbf{Claim 1:} $U, r^2 U_r$ are even in $\th$,  $r U_{\th}$   is  even  in $r$, and $\lim_{r\to 0}-r^2 U_r = rU$. In words, the last identity implies that the function $U$ is homogeneous of degree $-1$ on $r$  where $r$ is small. \\

The dynamics on $\mc H$ is thus 
\[ r'= \th_* \nu r, \ \nu'=\th_*(\frac{1}{2} \nu^2 + r^2 U_r)  \]
Since $u=0$, it is just the homothetic orbits  considered in Subsection \ref{subsec:homothetic}, but is regularized.  
There is one equilibrium point,  the intersection of the collision manifold and $\mc H$, denote it by $P$. The exact coordinates is 
\[   P=(0, -\nu_0, 0, 0), \frac{1}{2} \nu_0^2 = r U|_{r=0, \th=0}. \]

\subsection{The equilibrium point $P$ is hyperbolic } \label{subsec:P}

As in the Newtonian collinear three-body problem, the equilibrium $P$ is found to be hyperbolic \cite{Mcgehee1974}. 

We use the coordinates $r, u, \gamma$, and  the variable $\nu$ is  treated as a function of $r,u, \gamma$. The energy relation gives 
\begin{align*}
& \frac{\partial \nu}{\partial r}= \frac{(rU)_r-1}{\nu}= \frac{1}{\nu_0}, \  \frac{\partial \nu}{\partial u}=\th_* \frac{rU_\th \cos u -2\gamma^2 \sin u/\cos^5 u}{\nu}= 0, \ \\
& \frac{\partial \nu}{\partial \gamma}= \frac{-\gamma }{\nu \cos^4 u }= 0, 
\end{align*}
at the point $P$. 
 Then one finds that the linearized differential equations at $P$ have matrix 
 \begin{align*}
&\begin{bmatrix}
\th_*\nu  \cos^2 u + \th_*  r \cos^2 u \frac{\partial \nu}{\partial r}& \th_*  r  \cos^2 u \frac{\partial \nu}{\partial u} & \th_*  r  \cos^2 u \frac{\partial \nu}{\partial \gamma} \\
0&\gamma(\frac{1}{\cos u})'&\frac{1}{\cos u}\\
-\frac{\th_*}{2} \gamma \cos^2 u \frac{\partial \nu}{\partial r} + \th_* (r U_\th)_r\cos^4 u & \bigstar&  \frac{-2\sin u}{\cos^2 u}
\end{bmatrix}\\
=& \begin{bmatrix}
-\th_*\nu_0 & 0 & 0\\
0&0&1\\
0& \th_*^2 r U_{\th\th} & \frac{1}{2}\th_* \nu_0
\end{bmatrix}
 \end{align*}
 where $$\bigstar=-\frac{\th_*}{2} \gamma (\nu \cos^2 u)_u + rU_{\th\th}\th_*^2 \cos^5 u +  rU_{\th}\th_* (\cos^4 u)_u-2\gamma^2 (\frac{\sin u}{\cos^2 u})_u,$$ 
and we use the fact in Claim 1 and the following Claim 2. 

\textbf{Claim 2:} At the point $ r=0, u=0$, we have  $ r U_{\th\th} >0$.  \\

Thus, the equilibrium $P$ is hyperbolic, with eigenvalues $\lambda_1=-\th_* \nu_0<0$ and $\lambda_2<0, \lambda_3>0$. Then it has  two-dimensional stable manifold and  one-dimensional unstable manifold.  
The  eigenvectors are $(1,0,0)$, $(0, 1, \lambda_2)$, and $ (0, 1, \lambda_3)$. The first stable eigenvector is tangent to the  homothetic orbit $\mc H$.  Note that the other stable eigenvector$ (0, 1, \lambda_2)$ points out of $\mc W$ since $\gamma\ge 0$ in $\mc W$. 
 It follows that $\mc H\bigcap \mc W= W^s(P) \bigcap \mc W$. The unstable manifold of $P$ is on the collision manifold, with one branch in $\mc W$. 

\begin{lemma}\label{lem:branch}
The branch of $W^u(P)$ in $\mc W$ exits $\mc W$ at a point of the form  $(0, \nu, \frac{\pi}{2}, \gamma)$ with $\nu < 0$.
\end{lemma}

The following fact will be used in the proof. 

\textbf{Claim 3:} Restricted on $r=0$, the maximum of  $2rU \cos ^2 u$ is at $u=0$. \\

\begin{proof}
	Consider the system  for $u, \nu$.  By the energy relation, 
the equations read 
\begin{align*}
\nu'&=\th_*\cos^2 u(-\frac{1}{2} \nu^2 + 2 rh + 2rU + r^2 U_r)=\th_*\cos^2 u(rU -\frac{1}{2} \nu^2) \\
u '&=\frac{ \gamma}{\cos u}= \cos u \sqrt{  2 rU- \nu^2 }. 
\end{align*}
Then 
\[  \frac{d\nu}{du} =\frac{\th_*}{2} \sqrt {2rU \cos ^2 u-\nu^2 \cos^2 u}  \le  \frac{\th_*}{2} \sqrt {2rU \cos ^2 u}.  \]
So 
\[  \frac{d\nu}{du} \le  \frac{\th_*}{2} \nu_0 \]
which implies  hat the increment  in $\nu $ for $0 \le u \le \frac{\pi}{2}$  satisfies: 
\[\Delta \nu \le \frac{\pi}{2}  \frac{\th_*}{2} \nu_0   \le \nu_0\]
since $\th_* \le \frac{\pi}{4}$. Since the branch of $W^u(P)$ begins near $P$, and  $P$ has coordinates $u=0$ and $\nu=-\nu_0$,
then  it arrives at $u = \frac{\pi}{2}$ without  crossing the line  $\nu= 0$.
\end{proof}

\subsection{$\mc W$ is  a Wazewski set}

In this subsection, we identify  the subsets $\mc W_0, \mc E$ and show

\begin{theorem}\label{thm:Waze_set}
$\mc W$ is  a Wazewski set for the flow on the energy manifold.
\end{theorem}

The first property obviously holds since the set $\mc W$ is closed. For the second property, we first identify the subsets $\mc W_0, \mc E$.

\begin{lemma} \label{lem:leaving_set} 
$\mc W_0 =\{x\in \mc W: \exists t>0, \phi_t(x)\notin \mc W\}=\mc W\setminus \mc H$
\end{lemma}
The following facts will be used in the proof. 

\textbf{Claim 4:} $2rU\cos^2 u|_{u=\frac{\pi}{2}}$ has a positive lower bound $c_2^2$. 

\textbf{Claim 5:}  The function $\th_* r U_\th \frac{\cos^4 u}{\sin u} $ has a positive lower bound $c_3$. \\

\begin{proof}
	Let $x_0=(r_0, \nu_0, u_0, \gamma_0)\in \mc W$, It  is easy to that the solution begin from $x_0$ exist as long as it remains in $\mc W$. 
	Now suppose $x_0 \in \mc W \setminus  \mc H$. Our goal is to show  that $\phi_t(x_0)$ eventually leaves $\mc W$. If $u_0=0$ then $u'(0)=\gamma_0 >0$ since $x_0\notin \mc H$. It follows that for every $t_0 >0, u(t_0)>0.$  Thus it is enough to assume   $u_0 > 0.$
	
	Let  $u_0 $ be a positive constant and  $\mc W_{u_0} =\{x \in \mc W: u\ge u_0\}$. Since $u(t)$ is
	non-decreasing in $\mc W$, $\mc W_{u_0} $ is positively invariant relative to $\mc W$. We show  below that there are two  constants $c_0 > 0$ and $c_1 > 0$ such that for
	every  $x\in  \mc W_{u_0}$ 
	\[  {\it either }  \ \frac{\gamma}{\cos u }\ge c_0, \  {\it or}\ (\frac{\gamma}{\cos u })'\ge c_1. \]
	
	Then it is easy to see that 
	$\phi_t(x_0)$ must eventually leave $\mc W$.   Note that 
	\[ \mc W_{u_0} = \mc W_{u_0}^+\bigcup \mc W_{u_0}^-,     \]
	  where 
	  \[  W_{u_0}^+=\{  x\in W_{u_0},  \frac{\gamma}{\cos u }\ge c_0 \}, \  \mc W_{u_0}^- =\{  x\in W_{u_0},  0\le \frac{\gamma}{\cos u }<  c_0 \}.  \]
	Since $(\frac{\gamma}{\cos u })'\ge c_1>0$ in  $\mc W_{u_0}^-$
	 it implies that  an orbit segment can stay in $\mc W_{u_0}^-$ for time at most $c_0/c_1$, and then would enter $\mc W_{u_0}^+$.    Note that  $\mc W_{u_0}^+$ is positively invariant relative to $\mc W_{u_0}$. Finally, 
	an orbit can remain in $\mc W_{u_0}^+$ for  time not longer than $\frac{\pi}{2 c_0}$ since $u'=\frac{\gamma}{\cos u }\ge c_0$. Hence,  every  
	orbit starting  in $\mc W_{u_0}$ must  leaves $\mc W$ eventually.

We now  construct $c_0 >0, c_1 >0$ such that either $ \frac{\gamma}{\cos u }\ge c_0$ or $ (\frac{\gamma}{\cos u })'\ge c_1$ for all
	$x\in \mc W_{u_0}$.  For $u=\frac{\pi}{2}$, the equation \eqref{equ:energy_mcgehee_2} implies  $\frac{\gamma}{\cos u } = \sqrt{2rU\cos^2 u}|_{u=\frac{\pi}{2}}\ge c_2$.  We can  choose  $c_0$ to be  less than $c_2$ then  $ \frac{\gamma}{\cos u }\ge c_0$ holds for $u=\frac{\pi}{2}$.  For  $u_0\le u<\frac{\pi}{2}$,  we have 
	\begin{align*}
(\frac{\gamma}{\cos u })' &= \frac{\gamma'}{\cos u} + \tan u (\frac{\gamma}{\cos u })^2\\
&=-\frac{1}{2}\th_* \nu \gamma \cos u  + \th_* r U_\th \cos^3 u -2  \sin u \frac{\gamma^2}{\cos^3 u}+ \tan u (\frac{\gamma}{\cos u })^2\\
&\ge   \tan u \left( \th_* r U_\th \frac{\cos^4 u}{\sin u} - (\frac{\gamma}{\cos u })^2\right)\\
&\ge  \tan u_0 \left( c_3- (\frac{\gamma}{\cos u })^2\right). 
	\end{align*}
	Then we take $c_0$ such that $c_0 \le c_2, c_0^2\le  \frac{c_3}{2}$, and take  $c_1 = \frac{c_3 \tan u_0}{2}$. Then on $u=\frac{\pi}{2}$, we have $\frac{\gamma}{\cos u } \ge c_2 \ge c_0$. For  $u_0\le u<\frac{\pi}{2}$, if $\frac{\gamma}{\cos u } \le  c_0$, then 
	$(\frac{\gamma}{\cos u })' \ge   \tan u_0 \left( c_3- (\frac{\gamma}{\cos u })^2\right)\ge c_1$, as required. 
	
	\end{proof}

It remains to identify the immediate exit set $\mc E$.  It is useful  to distinguish two subsets of the boundary. Let $x = (r,\nu,u,\gamma)$ and let
\begin{align*}
\mc B_1 &= \{ x \in \mc W :  u = \frac{\pi}{2} \}, \\
\mc B_2 &=\{x\in \mc W: \nu=0, 0\le u <\frac{\pi}{2},   2 rU  + r^2 U_r -2r\ge 0\}.
\end{align*}
Obviously, the two subsets $\mc B_1$ and $\mc B_2$ are relatively closed in $\mc W$.  Hence Theorem \ref{thm:Waze_set} is proved once we show

\begin{lemma}
	The immediate exit set of $\mc W$ is $\mc E = \mc B_1 \bigcup \mc B_2$.
\end{lemma}

\begin{figure}[h!]
	\centering
	\includegraphics[scale=0.4]{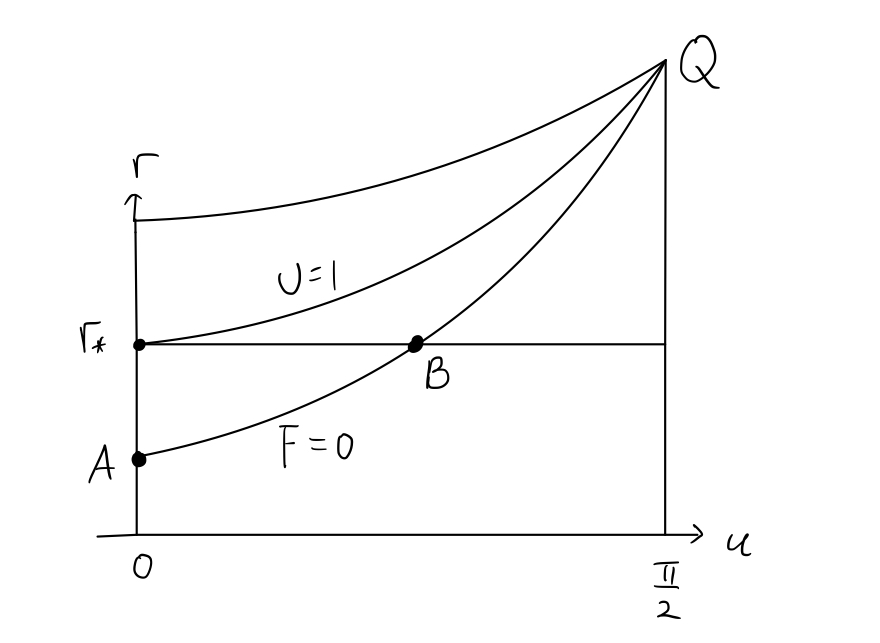}
	
	\caption{The configuration space and the rectangle $[0,\frac{\pi}{2}]\times [0,r_*]$ }
	\label{fig:rectangle}
\end{figure}
The following fact will be used in the proof. 

\textbf{Claim 6:}  Let $F(r, u)= 2rU + r^2 U_r - 2 r$.      $F_u>0$ for $0\le r\le r_*, 0<u\le \frac{\pi}{2}$. At $u=0, F=0$, we have $F_u=0, F_{uu}>0$.  
\\

\begin{proof}
	As claimed,  $\sqrt{2rU\cos^2 u}|_{u=\frac{\pi}{2}}$ has a positive lower bound, so $\mc B_1 \subset \mc E$.  Consider a point  $x\in \mc B_2$. Note that $h=-1, \nu=0$, then 
	\begin{align*}
	\nu'&=\th_*( \frac{\gamma^2}{\cos^2 u} +  r^2 U_r \cos^2 u )=\th_*\cos^2 u(  2rU + r^2 U_r - 2 r )\ge 0. 
	\end{align*}
	Let $F(r, u)= 2rU + r^2 U_r - 2 r$.  Since  $F_u>0$ in the rectangle   $(0,\frac{\pi}{2}]\times [0,r_*]$,  the set $F=0$ in this rectangle  is a curve  bounded by the two points $A,B$ (see Figure \ref{fig:rectangle}).

	The curve divides the rectangle into two parts. The bottom $r=0$ is in the set $F>0$ since on which $F=rU$, while the vertex $u=0, r=r_*$ is in the set $F<0$ since $U=1$ and that $ r^2 U_r<0$.

	If $F>0$, then $\nu'>0$ and $x$ is an immediate exit point. If $F=0$ and $u\ne 0$, one has $\nu=\nu'=0$ and one finds that the second derivative reduces to 
	\[  \nu''=\th_* (  -\sin 2 u F +\cos^2 u F_u) u' + \cos^2 u F_r r' = \cos u F_u  \gamma>0,  \]
	and $x$ is  an immediate exit point. Finally, if $u=\nu=0$ and $ F=0$, one has $\nu=\nu'=\nu''=0$, The third derivative at the point $A$  is found to be 
	\[  \nu'''=\th_* \gamma^2 F_{uu}>0.  \] 
	Again, 	$x$ is  an immediate exit point.

It remains to check  that there are no other immediate
	exit points.  Suppose that $x_0\in \mc W$ is an immediate exit point and it  is not in $\mc B_1 \bigcup \mc B_2$. Following the argument in \cite{Moeckel-Schubart}, it is enough the check the following cases. 
	
	First,  it may happen that $r_0=0$ but $r(t) < 0$ for small positive times. This is impossible because the manifold  $\{r = 0\}$ is invariant. 
	
	Secondly, it may happen  that $u_0=0$ but $u(t) < 0$ for small positive times. It  requires $u'(0)=\gamma_0 \le 0$ and since $x_0\in \mc W$ this means $\gamma_0 = 0$, so $x_0 \in \mc H$ and  points of $\mc H$ are certainly not leaving $\mc W$. 
	
Thirdly, it may happen  that $\nu_0=0$ but $\nu(t)$ increases for small positive times.
	This forces $\nu'(0)\ge 0$ and then $x_0 \in \mc B_2$.
	
	Fourthly, it may happen that $r_0=r_*$ but $r(t)$ increases.
	This forces $r'(0)= 0$ and then $\nu_0 =0$, i.e., the coordinates of the point is $r=r_*, 0\le u< \a$, where $(r_*, \a)$ is the coordinates of the point $B$, see Figure \ref{fig:rectangle}. So $\nu'<0$.  Then one finds 
\[ r''=\th_*( \nu (r \cos^2 u)' + r \cos^u \nu' )= \th_* r_* \nu'  \cos^2 u <0,  \]
	This mode of existing is impossible. 

At last, it may happen that  $\gamma_0=0$ but $\gamma(t)$ decreases for small positive times. If $u_0=0$, then $x_0\in \mc H$, and points of $\mc H$ are certainly not leaving $\mc W$. If $u_0=\frac{\pi}{2}$, then $x_0 \in \mc B_1$.  One may assume $0<u_0 <\frac{\pi}{2}$. In this case, it follows from the proof of Lemma  \ref{lem:leaving_set} that there are positive constants $c_0, c_1$ such that $(\frac{\gamma}{\cos u })'\ge c_1 > 0$ whenever $\frac{\gamma}{\cos u }<c_0$. So this mode of exiting is also impossible. This completes the proof.
	
	\end{proof}




\subsection{The shooting argument}

Finally, we can complete the construction  of the symmetric periodic orbit. Recall that it suffices to construct the first quarter, which is required to be  an orbit from $u=0, \nu=0$ to  $u=\frac{\pi}{2}, \nu=0$.  

Since  $\mc W$ is a Wazewski set,   the time required to reach $\mc E$ depends continuously on initial conditions and so there is a continuous flow-defined map from $\mc W_0$ to $\mc E$. 
 The map is also continuous if we restrict the domain to 
\[  \mc S=\{  (r, \nu, u, \gamma ) \in \mc W_0, u = \nu = 0, 0\le r< r_*\}. \]
That is, the flow-defined map $F:  \mc S\to  \mc E$ is continuous. 
Let 
\[ \mc T=\{   (r, \nu, u, \gamma )  \in \mc W, u =\frac{\pi}{2},  \nu = 0\}. \]
Note that $\mc T\subset \mc E$ and  that  $\mc S$ and $\mc T$ are two of the edges in the boundary of the three-dimensional Wazewski set $\mc W$ (shown as bold vertical lines in Figure \ref{fig:Waze-set}). Then the construction of the first quarter of the orbit reduces to show that 
 \[   F(\mc S) \bigcap  \mc T\ne  \emptyset.\]

First, note that  part of $\mc S$ near $r=0$ is contained in $\mc B_2\subset \mc E$. These points exit $\mc W$ immediately, so  the map $F$ is the identity there. Secondly,  points of $\mc S$ with $r$ close to $r_*$ will enter the interior of $\mc W$ and exit elsewhere. By continuous dependence of the initial conditions,  these points will follow the  homothetic orbit $\mc H$ to a neighborhood of the  equilibrium point  $P = (0, -\nu_0, 0, 0)$.  Then the lambda lemma \cite{Wiggins1990} implies that  they will follow a branch of the unstable manifold $W^u(P)$, which is one-dimensional and is contained entirely in the invariant manifold $r=0$, as shown in Subsection \ref{subsec:P}. 
Furthermore, by Lemma   \ref{lem:branch},  one of the two branches lies  in $\mc W_0$ and it  goes to some point on $r=0, \nu <0, u=\frac{\pi}{2}$.  Then the lambda lemma implies that the image of  points near the upper endpoint of $\mc S$ under the continuous mapping $F$ are on $\mc B_1\setminus \mc T$.

We can now complete the shooting argument. Recall that there is continuous map $F: \mc S \to \mc E$, $\mc E=\mc B_1 \bigcup \mc B_2$ and that  $\mc B_1$and $\mc B_2$  are two-dimensional continuum  meeting along the edge $\mc T$.   As we have shown,  the image of points near $r=0$ under  $F$ are in $\mc B_2\setminus \mc T$, while the image of points near $r=r_*$ under  $F$ are  in $\mc B_1\setminus \mc T$.
 It follows that there must exist at least one intersection point $U\in F(\mc S) \bigcap \mc T$. This shows that $ F(\mc S) \bigcap  \mc T\ne  \emptyset$ and  completes the existence proof for the symmetric periodic orbits.

\begin{remark}
	The orbit constructed lies in the energy manifold $h=-1$. By restricting the configuration to the rectangle  $0\le r\le r_*, 0\le u\le \frac{\pi}{2}$, where  $u=0, r=r_*$ is the intersection of $U=1$ and $u=0$, the six Claims hold.  Consider an energy manifold $h< -1$. Since $U$ is decreasing on $u=0$, the  intersection of $U=-h$ and $u=0$ is lower than $r_*$. Then  the six Claims made in this section still hold and all arguments can be applied as well.  Thus, we have the following 
\end{remark}

\begin{theorem}\label{thm:main}
 Given three positive masses $m_1 = m_2$ and $m_3$ and an  energy $h\le -1$. 
 Then there exists a symmetric periodic solution of the collinear three-body problem on $\S^1$ with energy $h$ and regularized double collisions. The orbit has the following features.
 \begin{itemize}
 	\item The configuration lies in region I.
 	\item   In  the first quarter of the orbit, the masses move from the Eulerian central configuration with $m_3$ in the middle of  $m_1, m_2$ to a double collision between $m_2$ and $ m_3$. At the moment of the double collision  the velocity of $m_1$ is zero.
 	\item The second quarter of the orbit is the time-reverse of the first, and the second half is the reflection of the first half with the roles of $m_1$ and  $m_2$ reversed.
 \end{itemize} 
 
\end{theorem}

 \section{Appendix: Proofs of the six Claims}

 Recall that     
 \begin{align*}
 rU&= n^2 r\cot (rA_1 \cos \th )  +m n \left(r\cot[ r  A_2 \sin(\th+\th_*)] +   r\cot[ r  A_2 \sin(\th_*-\th)] \right)\\
 rU_\th&=n^2\frac{r^2 A_1 \sin \th}{\sin^2 (rA_1 \cos \th )}  +m n \left(  -\frac{r^2  A_2 \cos (\th + \th_*)}{\sin^2[ r  A_2 \sin(\th+\th_*)]} +    \frac{r^2  A_2 \cos (\th_*- \th)}{\sin^2[ r  A_2 \sin(\th_*-\th)]} \right)\\
 r^2 U_r&=-n^2\frac{r^2 A_1 \cos \th}{\sin^2 (rA_1 \cos \th )}   -m n \left(  \frac{r^2  A_2 \sin (\th + \th_*)}{\sin^2[ r  A_2 \sin(\th+\th_*)]} +    \frac{r^2  A_2 \sin (\th_*- \th)}{\sin^2[ r  A_2 \sin(\th_*-\th)]} \right)
 \end{align*}

\subsection{}
 \textbf{Claim 1:} $U, r^2 U_r$ are even in $\th$,  $r U_{\th}$   is  even  in $r$, and $\lim_{r\to 0}-r^2 U_r = rU$. In words, the last identity implies that the function $U$ is almost homogeneous of degree $-1$ on $r$  where $r$ is small. 
 
It is easy to see by the explicit form of the functions.

  \subsection{}   \textbf{Claim 2:} At the point $ r=0, u=0$, we have  $ r U_{\th\th} >0$.

We show that      
     \[  rU_\th=n^2\frac{r^2 A_1 \sin \th}{\sin^2 (rA_1 \cos \th )}  +m n \left(  -\frac{r^2  A_2 \cos (\th + \th_*)}{\sin^2[ r  A_2 \sin(\th+\th_*)]} +    \frac{r^2  A_2 \cos (\th_*- \th)}{\sin^2[ r  A_2 \sin(\th_*-\th)]} \right) \]
     is strictly increasing on $\th$ at $ r=0, \th=0$. 
     
    The first term is strictly increasing in $\th$, at $r=\th=0$. Direct  computation gives 
     \begin{align*}
 (\frac{r^2 A_1 \sin \th}{\sin^2 (rA_1 \cos \th )}  )_\th&= r^2 A_1\frac{ \sin^2 (rA_1 \cos \th )\cos \th -  r A_1 \sin^2 \th \sin ( 2 rA_1 \cos \th ) }{\sin^4 (rA_1 \cos \th )}\\
 &\to \frac{1}{A_1}>0.   
    \end{align*}

The second term, denoted by $g(r,\th)$,  is  an   increasing function  on $\th$ in a neighborhood of $r=\th=0$. Indeed,  
 $g(r, 0)=0$, and $g(r, \th)>0$ if $0< \th\le \th_*$ and $r$ is small. Note that $\th_*\le \frac{\pi}{4}$, then 
 \[\cos (\th_*- \th) >\cos (\th_*+ \th), \  r A_2\sin (\th_*- \th) < r A_2\sin (\th_*+ \th)<\frac{\pi}{2}.  \] 
 Thus, 
 \[   g(r, \th)=  m n r^2  A_2  \left(  \frac{ \cos (\th_*- \th)}{\sin^2[ r  A_2 \sin(\th_*-\th)]} -\frac{\cos (\th + \th_*)}{\sin^2[ r  A_2 \sin(\th+\th_*)]} \right)>0,    \]  
and the derivative $g_\th(0,0)\ge 0$.

      \subsection{} \textbf{Claim 3:} Restricted on $r=0$, the maximum of  $2rU \cos ^2 u$ is at $u=0$. 

Recall that 
\[  rU= n^2 r\cot (rA_1 \cos \th )  +m n \left(r\cot[ r  A_2 \sin(\th+\th_*)] +   r\cot[ r  A_2 \sin(\th_*-\th)] \right)\]
Let $r\to 0$, we have 
\begin{align*}
rU \cos ^2 u &=  n^2 \frac{\cos^2 u}{ A_1 \cos(\th_* \sin u)}   +m n  \cos^2 u\left(\frac{1}{ A_2 \sin (\th_*  + \th)}+\frac{1}{ A_2 \sin (\th_*  - \th) }\right)\\
&=n^2 \frac{\cos^2 u}{ A_1 \cos(\th_* \sin u)}   +\frac{2mn \sin \th_*}{A_2} \frac{ \cos^2 u\cos \th }{ \cos^2 \th - \cos^2 \th_* }. 
\end{align*}

The first term is a decreasing function of $u$. Since  $\th_* \le \frac{\pi}{4}<1$, we have 
  \begin{align*}
(\frac{\cos^2 u}{  \cos(\th_* \sin u)}  )'&= \frac{\cos u}{\cos^2 (\th_*\sin u) } [\th_*\cos ^2u \sin (\th_*\sin u)-2 \sin u \cos (\th_*\sin u) ]\\
&\le \frac{\cos u}{\cos^2 (\th_*\sin u) } 2[\cos u  \sin (\th_*\sin u)-\sin u \cos (\th_*\sin u) ]\\
& =  \frac{\cos u}{\cos^2 (\th_*\sin u) } 2[ \sin ( \th_*\sin u-u )] \le 0. 
 \end{align*}
	
It remains to show 
	\begin{equation}\label{equ:the_inequality}
	\frac{1 }{ 1 - \cos^2 \th_* }\ge \frac{ \cos^2 u\cos \th }{ \cos^2 \th - \cos^2 \th_* }, \  u\in [0,\frac{\pi}{2}]. 
	\end{equation}
it is equivalent to 
\[  \mc J=\cos^2 \th - \cos^2 \th_*  -(1 - \cos^2 \th_*) \cos^2 u\cos \th\ge 0.  \]	
	View $\mc J$ as a function of the two variables $(\th_*, \th)$, on the triangular region  $0< \th_*\le  \frac{\pi}{4}, 0\le \th\le \th_*$. Note that $\mc J(\th_*, \th_*)=0$ and 
	\[ \frac{\partial \mc J}{\partial \th_*} = 2 \sin \th_* \cos \th_* (1-\cos^2 u\cos \th)\ge 0. \]	
	We conclude that the function $\mc J$ is non-negative on the triangular region.


	\subsection{} \textbf{Claim 4:} $2rU\cos^2 u|_{u=\frac{\pi}{2}}$ has a positive lower bound. 
	
       Recall that $rU = n^2 r \cot d_{12} + mn (\cot d_{13} + \cot d_{23})$, and the fact that we are not at the singularity $Q$.  Hence, when $u=\frac{\pi}{2}$,  we have $d_{23}=0$, and  the two distance $d_{12}, d_{13}$ are different from $0, \pi$. 

Hence 
\[  rU\cos^2 u = mn \cot d_{23} \cos^2 u =  mn r \cot [r A_2 \sin (\th_*-\th)] \cos^2 u  \]
Since $\cos^2 u= \sin^2 (\frac{\pi}{2}-u)$, and $\th_*-\th=2\th_* \sin^2 (\frac{\pi}{2}-u)$, so we obtain 
\[rU\cos^2 u  =   \frac{mn}{2\th_*A_2 }.   \]

   \subsection{}  \textbf{Claim 5:}  The function $  r U_\th \frac{\cos^4 u}{\sin u} $ has a positive lower bound. 
 
 Recall that      
 \[  rU_\th=n^2\frac{r^2 A_1 \sin \th}{\sin^2 (rA_1 \cos \th )}  +m n \left(  -\frac{r^2  A_2 \cos (\th + \th_*)}{\sin^2[ r  A_2 \sin(\th+\th_*)]} +    \frac{r^2  A_2 \cos (\th_*- \th)}{\sin^2[ r  A_2 \sin(\th_*-\th)]} \right). \]
 
 We first claim that the function $rU_\th$  is non-negative. The first term is non-negative. For the second term, which has been denoted by $g(r,\th)$.   We have showed that $g(r, 0)=0$, and $g(r, \th)>0$ if $0< \th\le \th_*$ and $r$ is small. Now we show that 
 \[   0\le r\le r_*, \ 0< \th\le \th_*, \Rightarrow g(r, \th)\ge 0. \]
For this, it suffices to show that 
\[  r_* A_2 \sin (2\th_*)\le \frac{\pi}{2}.  \]
 
 Recall that at $\th=0, r=r_*$, we have $U=1$.  Note that $r A_2 \sin \th_* = r \sqrt{\frac{m+1}{2nm}} \sqrt{\frac{m}{m+1}}= rA_1/2$. Then 
 \[   n^2 \cot (r_*A_1)  + 2 mn \cot (r_* A_1/2)=1.  \]
 Let $a=\cot (r_* A_1/2)$.  Note that  $a>0$ since $r A_2 \sin \th_*=d_{23}<\frac{\pi}{2}$.  Then
 \[  2m n a +n^2 \frac{a^2-1}{2a} =1, \Rightarrow  (4m n +n^2)  a^2-2a-n^2=0.   \]
 So 
 \[ a=\cot (r_* A_1/2)=\frac{4}{-7 m^2+6 m+1}+\frac{\sqrt{-7 m^4+20 m^3-18 m^2+4 m+17}}{-7 m^2+6 m+1}.   \]
Let $g(m)=-7 m^4+20 m^3-18 m^2+4 m+17$.  We have $g'=-4 (m-1)^2 (7 m-1)$, so $g(m)\ge \min \{g(0), g(1) \}=16$.  Since $-7 m^2+6 m+1\le \frac{16}{7}$, we obtain the desired estimate 
\begin{align}
&\label{equ:estimate}\cot (r_* A_1/2) \ge \frac{7}{2}, \Rightarrow    r_* A_1/2 \le  \frac{\pi}{10}\\
&r_* A_2 \sin (2\th_*)= r_* A_1\cos \th_*= \frac{r_* A_1}{\sqrt{m+1}}<\frac{\pi}{5}.  \notag
\end{align}

 Now we show that $\th_* r U_\th \frac{\cos^4 u}{\sin u}, 0<u_0<u<\frac{\pi}{2} , 0\le r\le r_*$, has a positive lower bound.  
 Obviously,  at $u=\frac{\pi}{2}$, the second term  equals to 
 \begin{align*}
& \lim_{ u\to \frac{\pi}{2}} \frac{r^2  A_2 \cos (\th_*- \th)}{\sin^2[ r  A_2 \sin(\th_*-\th)]}  \frac{\cos^4 u}{\sin u}  \\
&=\lim_{ u\to \frac{\pi}{2}} \frac{mn  \sin^4(\pi/2-u)  }{ A_2 \sin^2[ 2\th_* \sin^2(\frac{\pi/2-u}{2})]} =\frac{4mn}{A_2 \th_*^2}. 
 \end{align*}
 Then there is some $u_1<\frac{\pi}{2}$ such that 
 \[  \th_* r U_\th \frac{\cos^4 u}{\sin u} \ge  \frac{2mn}{A_2 \th_*^2}, \ u_1\le u\le \frac{\pi}{2}, 0\le r\le r_*. \]
 For the first term, let $\th_0= u_0 \sin u$, then 
 \[ \frac{r^2 A_1 \sin \th}{\sin^2 (rA_1 \cos \th )}  \frac{\cos^4 u}{\sin u}  \ge  \frac{r^2 A_1^2 \sin \th \cos^4 u}{ A_1\sin^2 (rA_1 )}  \ge \frac{ \sin \th_0 \cos^4 u_1}{A_1}, \   \ u_0\le u\le u_1, 0\le r\le r_*. \]
 Thus, we conclude  that $ r U_\th \frac{\cos^4 u}{\sin u}, 0<u_0<u<\frac{\pi}{2} , 0\le r\le r_*$, has a positive lower bound.  
 
 \subsection{}  \textbf{Claim 6:} Let $F(r, u)= 2rU + r^2 U_r - 2 r$.      $F_u>0$ for $0\le r\le r_*, 0<u\le \frac{\pi}{2}$. At $u=0, F=0$, we have $F_u=0, F_{uu}>0$.  
  
  Recall that 
  \begin{align*}
  2rU +r^2 U_r&= n^2 r [ 2\cot (rA_1 \cos \th ) -\frac{r A_1 \cos \th}{\sin^2 (rA_1 \cos \th )}  ]  +m n r \{ 2 \cot[ r  A_2 \sin(\th_*-\th)]\\
  &-   \frac{r  A_2 \sin (\th_*- \th)}{\sin^2[ r  A_2 \sin(\th_*-\th)]} +2\cot[ r  A_2 \sin(\th+\th_*)]  -\frac{r  A_2 \sin (\th + \th_*)}{\sin^2[ r  A_2 \sin(\th+\th_*)]}  \}
  \end{align*}
  Introduce  new variables
  \[ \rho= rA_1 \cos \th, \  \xi = rA_2 \sin (\th_* -\th), \  \eta= rA_2 \sin (\th_* +\th),     \]
and define $f(x)= 2\cot x -\frac{x}{\sin^2 x}$.   Then 
  \[ 2rU +r^2 U_r= n^2 r f(\rho) + mn r [f(\xi) + f(\eta)].   \]
  
 Let us first study the function $f(x)$.  On $[0, \frac{\pi}{5}]$
 \begin{align*}
f'(x)&= -\frac{3-2 x \cot x }{\sin ^2 x} <0, \\
f''(x)&= -\frac{2}{\sin^4 x} (2 x-2 \sin  2 x+x \cos 2 x ) \ge   \frac{2}{\sin^4 x} (2 \sin  2 x-3x )\ge 0,
 \end{align*}
since  $k(x)=2 \sin  2x  -3x$ is  a concave function on $[0, \frac{\pi}{2}]$.  Thus its value on $[0, \frac{\pi}{5}]$ is at least $\min\{ k(0),k(\frac{\pi}{5})  \}=0$. 
\begin{align*}
&\rho'=-rA_1 \sin \th,    &\xi'=-rA_2 \cos(\th_*-\th), \ \  &\eta'=rA_2 \cos(\th_*+\th), \\
 &\rho''=-\rho,  &\xi''=-\xi, \ \ \ \ \ \ \ \  &\eta''=-\eta, \\
  & \rho\le \frac{\pi}{5}, &\xi \le \eta \le r_* A_2 \sin 2 \th_* < \frac{\pi}{5}. 
 \end{align*}

  For the first derivative, one finds $F_{u} =  F_\th\th_* \cos u$, then it suffices to show that $F_\th = (2rU + r^2 U_r )_\th >0$ for $0\le r\le r_*, 0<\th< \th_*$. 
  \[   F_\th =   n^2 rf'(\rho) \rho' +mn r (f'(\xi) \xi'  +f'(\eta) \eta'  ).   \] 
The first term is positive  if $\th \in (0, \frac{\pi}{5}]$, and it is zero if $\th=0$. The second term is  zero if $\th=0$, and it is positive  if $\th \in (0, \frac{\pi}{5})$ since both $-f'(x)$ and $\cos(x)$ are decreasing.  Hence, we have proved that  $F_u>0$ for $0\le r\le r_*, 0<u< \frac{\pi}{2}$ and  $F_u=0$ for $0\le r\le r_*, u=0$. 

For the second derivative at the point $A$, one finds
\[   F_{uu} =  F_{\th\th} (\th_* \cos u)^2 - F_\th \th_* \sin u=  F_{\th\th} \th_* ^2.    \]
and
\begin{align*}
F_{\th\th}&= n^2 r[f''(\rho) (\rho')^2 - f'(\rho) \rho ]+mn r [f''(\xi) (\xi')^2 - f'(\xi) \xi + f''(\eta) (\eta')^2 - f'(\eta)\eta  ]\\
&=n^2 r[- f'(\rho) \rho ]+mn r [ 2f''(\xi) (\xi')^2 - 2f'(\xi) \xi  ]>0. 
\end{align*}
Hence, we have proved that  $F_{uu}>0$ at the point $A$.

\textbf{Acknowledgments. }  The author would like to thank Cristina Stoica and Jean-Marie Becker for enlightening discussions.


\begin{thebibliography}{10}
	
	\bibitem{borisov2018reduction}
	AV~Borisov, LC~Garc{\'\i}a-Naranjo, IS~Mamaev, and James Montaldi.
	\newblock Reduction and relative equilibria for the two-body problem on spaces
	of constant curvature.
	\newblock {\em Celestial Mechanics and Dynamical Astronomy}, 130(6):1--36,
	2018.
	
	\bibitem{Conley_isolating}
	C.~Conley and R.~Easton.
	\newblock Isolated invariant sets and isolating blocks.
	\newblock {\em Trans. Amer. Math. Soc.}, 158:35--61, 1971.
	
	\bibitem{Conely_retrograde}
	C.~C. Conley.
	\newblock The retrograde circular solutions of the restricted three-body
	problem via a submanifold convex to the flow.
	\newblock {\em SIAM J. Appl. Math.}, 16:620--625, 1968.
	
	\bibitem{Conly_book}
	Charles Conley.
	\newblock {\em Isolated invariant sets and the {M}orse index}, volume~38 of
	{\em CBMS Regional Conference Series in Mathematics}.
	\newblock American Mathematical Society, Providence, R.I., 1978.
	
	\bibitem{diacu2011}
	Florin Diacu.
	\newblock On the singularities of the curved n-body problem.
	\newblock {\em Transactions of the American Mathematical Society},
	363(4):2249--2264, 2011.
	
	\bibitem{diacumemoir}
	Florin Diacu.
	\newblock Relative equilibria in the 3-dimensional curved {$n$}-body problem.
	\newblock {\em Mem. Amer. Math. Soc.}, 228(1071):vi+80, 2014.
	
	\bibitem{diacu2012n2}
	Florin Diacu, Ernesto P{\'e}rez-Chavela, and Manuele Santoprete.
	\newblock The n-body problem in spaces of constant curvature. part ii:
	Singularities.
	\newblock {\em Journal of nonlinear science}, 22(2):267--275, 2012.
	
	\bibitem{diacu2018central}
	Florin Diacu, Cristina Stoica, and Shuqiang Zhu.
	\newblock Central configurations of the curved n-body problem.
	\newblock {\em Journal of Nonlinear Science}, 28(5):1999--2046, 2018.
	
	\bibitem{Easton1970}
	Robert~W. Easton.
	\newblock On the existence of invariant sets inside a submanifold convex to a
	flow.
	\newblock {\em J. Differential Equations}, 7:54--68, 1970.
	
	\bibitem{Mcgehee1974}
	Richard McGehee.
	\newblock Triple collision in the collinear three-body problem.
	\newblock {\em Invent. Math.}, 27:191--227, 1974.
	
	\bibitem{Moeckel-Schubart}
	Richard Moeckel.
	\newblock A topological existence proof for the {S}chubart orbits in the
	collinear three-body problem.
	\newblock {\em Discrete Contin. Dyn. Syst. Ser. B}, 10(2-3):609--620, 2008.
	
	\bibitem{Moser1973}
	J\"{u}rgen Moser.
	\newblock {\em Stable and random motions in dynamical systems}.
	\newblock Annals of Mathematics Studies, No. 77. Princeton University Press,
	Princeton, N. J.; University of Tokyo Press, Tokyo, 1973.
	\newblock With special emphasis on celestial mechanics, Hermann Weyl Lectures,
	the Institute for Advanced Study, Princeton, N. J.
	
	\bibitem{Schubart1956}
	J.~Schubart.
	\newblock Numerische {A}ufsuchung periodischer {L}\"{o}sungen im
	{D}reik\"{o}rperproblem.
	\newblock {\em Astronom. Nachr.}, 283:17--22, 1956.
	
	\bibitem{Shibayama2011}
	Mitsuru Shibayama.
	\newblock Minimizing periodic orbits with regularizable collisions in the
	{$n$}-body problem.
	\newblock {\em Arch. Ration. Mech. Anal.}, 199(3):821--841, 2011.
	
	\bibitem{Wazewski1947}
	Tadeusz Wa\.{z}ewski.
	\newblock Sur un principe topologique de l'examen de l'allure asymptotique des
	int\'{e}grales des \'{e}quations diff\'{e}rentielles ordinaires.
	\newblock {\em Ann. Soc. Polon. Math.}, 20:279--313 (1948), 1947.
	
	\bibitem{Wiggins1990}
	Stephen Wiggins.
	\newblock {\em Introduction to applied nonlinear dynamical systems and chaos},
	volume~2 of {\em Texts in Applied Mathematics}.
	\newblock Springer-Verlag, New York, 1990.
	
	\bibitem{Xia5body}
	Zhihong Xia.
	\newblock The existence of noncollision singularities in {N}ewtonian systems.
	\newblock {\em Ann. of Math. (2)}, 135(3):411--468, 1992.
	
\end{thebibliography}
\end{document}